\newcount\Review
\documentclass[conference]{IEEEtran}
\IEEEoverridecommandlockouts

\usepackage{amsmath, bm}
\usepackage{amssymb}
\usepackage{amsfonts}

\usepackage{graphicx}
\usepackage{cite}
\usepackage{picinpar}
\usepackage{url}
\usepackage{flushend}
\usepackage{colortbl}
\usepackage{soul}
\usepackage{multirow}
\usepackage{pifont}
\usepackage{color}
\usepackage{alltt}
\usepackage[hidelinks]{hyperref}
\usepackage{enumerate}
\usepackage{breakurl}

\usepackage[squaren]{SIunits}			
\usepackage{psfrag}							
\usepackage{upgreek}
\usepackage{placeins}
\usepackage{float}
\usepackage{printlen}

\newcommand{\R}{\mathbb{R}}
\newcommand{\Kcal}{\mathcal{K}}

\newcommand{\capname}{Sec. }
\newcommand{\figref}[1]{\figurename~\ref{#1}}
\newcommand{\tabref}[1]{\tablename~\ref{#1}}
\newcommand{\capref}[1]{\capname~\ref{#1}}

\newtheorem{remark}{Remark}

\begin{document}

\begin{center}
	\includegraphics[keepaspectratio=true, scale=1.08]{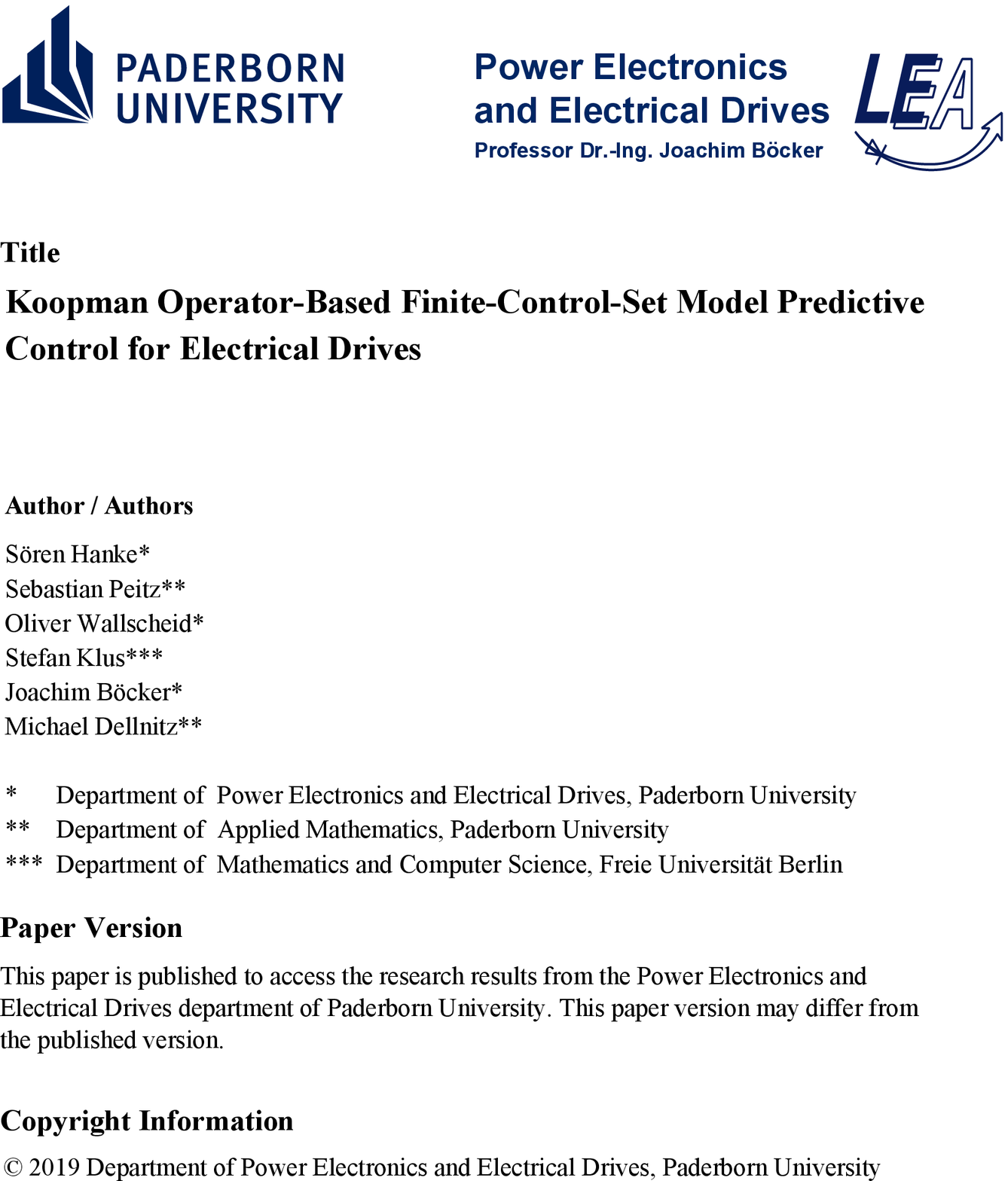}
\end{center}

\newpage

\title{Koopman Operator-Based Finite-Control-Set Model Predictive Control for Electrical Drives
\thanks{This work was funded by the German Research Foundation (DFG) under the reference number BO 2535/11-1.}
}

\author{

\IEEEauthorblockN{S{\"o}ren Hanke\IEEEauthorrefmark{1}, 
                  Sebastian Peitz\IEEEauthorrefmark{2}, 
									Oliver Wallscheid\IEEEauthorrefmark{1}, 
									Stefan Klus\IEEEauthorrefmark{3}
									Joachim B{\"o}cker\IEEEauthorrefmark{1} 
                  and Michael Dellnitz\IEEEauthorrefmark{2}}
													
\IEEEauthorblockA{\IEEEauthorrefmark{1}Power Electronics and Electrical Drives, Paderborn University, 33095 Paderborn, Germany,\\hanke@lea.upb.de, wallscheid@lea.upb.de, boecker@lea.upb.de}
\IEEEauthorblockA{\IEEEauthorrefmark{2}Applied Mathematics, Paderborn University, 33095 Paderborn, Germany,\\ speitz@math.upb.de, dellnitz@upb.de}
\IEEEauthorblockA{\IEEEauthorrefmark{3}Mathematics and Computer Science, Freie Universit\"at Berlin, 14195 Berlin, Germany,\\ stefan.klus@fu-berlin.de}
}

	%
%
	%
	%

\maketitle

\begin{abstract}
Predictive control of power electronic systems always requires a suitable model of the plant. Using typical physics-based white box models, a trade-off between model complexity (i.e. accuracy) and computational burden has to be made. This is a challenging task with a lot of constraints, since the model order is directly linked to the number of system states. Even though white-box models show suitable performance in most cases, parasitic real-world effects often cannot be modeled satisfactorily with an expedient computational load. Hence, a Koopman operator-based model reduction technique is presented which directly links the control action to the system's outputs in a black-box fashion. The Koopman operator is a linear but infinite-dimensional operator describing the dynamics of observables of nonlinear autonomous dynamical systems which can be nicely applied to the switching principle of power electronic devices. Following this data-driven approach, the model order and the number of system states are decoupled which allows us to consider more complex systems. Extensive experimental tests with an automotive-type permanent magnet synchronous motor fed by an IGBT 2-level inverter prove the feasibility of the proposed modeling technique in a finite-set model predictive control application.    
\end{abstract}


\definecolor{limegreen}{rgb}{0.2, 0.8, 0.2}
\definecolor{forestgreen}{rgb}{0.13, 0.55, 0.13}
\definecolor{greenhtml}{rgb}{0.0, 0.5, 0.0}

\section{Introduction}


Linear feedback control is the most frequently used control strategy employed in power electronic and drive applications because  of  its  simplicity and well-known design rules. However, those applications require to address nonlinear influences due to state and control action constraints, parameter changes (e.g. due to iron saturation) or in some cases even completely nonlinear control plants (e.g. LLC resonant converter). Hence, when designing linear feedback control loops, e.g. by utilizing classic PID elements, the control engineer has to manually tune the control parameters after a first analytical design step to ensure stability and acceptable performance. In addition, nonlinear control elements like anti-windup reset measures have to be added to the control loops. Consequently, most linear feedback control approaches are transformed step-wise into hand-tailored nonlinear control systems during the design process requiring highly experienced and application-specific engineering knowledge.   

In contrast, model predictive control (MPC) techniques inherently allow to address nonlinear model plants and manifold system constraints by suitable definition of an optimal control problem on a receding horizon. In general, MPC calculates the control output by minimizing a cost function that describes the desired system behavior. The cost function evaluates the model-based predicted system output with a reference trajectory. For each sampling instant, the MPC calculates a control action sequence that minimizes the cost function, yet only the first component of this sequence is applied to the system. 

MPC theory was intensively investigated throughout the last century and first applied to complex chemical processes where standard linear feedback control delivered only unsatisfactory results \cite{Morari1999} leading to first publications during the 1970s and 1980s \cite{Cutler1980}. First chemical processes considered for MPC had vast dominant time constants in the range of minutes or even hours and, therefore, sufficiently large time intervals were available to solve the optimization problems with the limited computational power of that time. Due to the increasing computational performance of digital signal processing units in the last three decades, MPC techniques also became feasible for electrical power systems with typical time constants in the range of milliseconds or even microseconds \cite{Cortes2008,Kouro2015,Vazquez2014} as well as for complex systems such as autonomous vehicles \cite{EPS+16,PSOB+17}.

For power electronic applications, one has to distinguish between finite-control-set (FCS-MPC) and continuous-control-set (CCS-MPC) approaches \cite{Vazquez2017}. In the latter, the actuating variable is an element of a continuous space, and the computed control signals have to be forwarded to a modulator (e.g. space-vector- or pulse-width-modulation) to receive the desired switching sequence which can then be applied to the semiconductor's driver circuits. The main advantages of CCS-MPC are a guaranteed fixed switching frequency and the usage of long prediction horizons since the controller turnaround-time is decoupled from the modulation carrier frequency. On the contrary, FCS-MPC directly computes the switching signal sequence suitable for driving power electronic devices. On the one hand, FCS-MPC gives additional degrees of freedom to the control problem, e.g. to find loss-optimal pulse patterns or to manipulate the frequency-spectrum behavior, but on the other hand, the controller turnaround-time has to be significantly smaller compared to CCS-MPC due to the missing modulation step. Furthermore, the underlying optimization problem is of combinatorial nature. Consequently, the FCS-MPC optimization routines have to be very fast and efficient which inherently requires streamlined plant models with minimal computational complexity. To achieve suitable control performance both high model accuracy and long prediction horizons (i.e. lightweight plant models) are required \cite{Geyer2015}. As these two model characteristics are obviously opposing goals, a trade-off decision has to be made.

\begin{figure}[ht]
	\centering
	\psfrag{CO}[cc][cc][1][0]{Control structure}
	\psfrag{FCS}[cc][cc][1][0]{FCS}
	\psfrag{MPC}[cc][cc][1][0]{MPC}
	\psfrag{T}[cc][cc][1][0]{$T,\omega$}
	\psfrag{PLL}[cc][cc][1][0]{PLL}
	\psfrag{rst}[cc][cc][1][0]{$i_{\mathrm{rst},i}$}
	\psfrag{u}[cc][cc][1][0]{$u_{\mathrm{dc},i}$}
	\psfrag{PA}[cc][cc][1][0]{Park}
	\psfrag{OS}[cc][cc][1][0]{OPS}
	\psfrag{P}[cc][cc][0.95][0]{IPMSM}
	\psfrag{d}[cc][cc][1][0]{$i_{\mathrm{d},i}$}
	\psfrag{q}[cc][cc][1][0]{$i_{\mathrm{q},i}$}
	\psfrag{dw}[cc][cc][1][0]{$i_{\mathrm{d},i}^*$}
	\psfrag{qw}[cc][cc][1][0]{$i_{\mathrm{q},i}^*$}
	\psfrag{sa}[cc][cc][1][0]{$s_{\mathrm{a},i}$}
	\psfrag{sb}[cc][cc][1][0]{$s_{\mathrm{b},i}$}
	\psfrag{sc}[cc][cc][1][0]{$s_{\mathrm{c},i}$}
	\psfrag{Tw}[cc][cc][1][0]{$T^*$}
	\psfrag{w}[cc][cc][1][0]{$\omega_{\mathrm{el},i}$}
	\psfrag{eps}[cc][cc][1][0]{$\varepsilon_{\mathrm{el},i}$}
	\ifnum\Review=0
		\includegraphics[width=0.9\columnwidth]{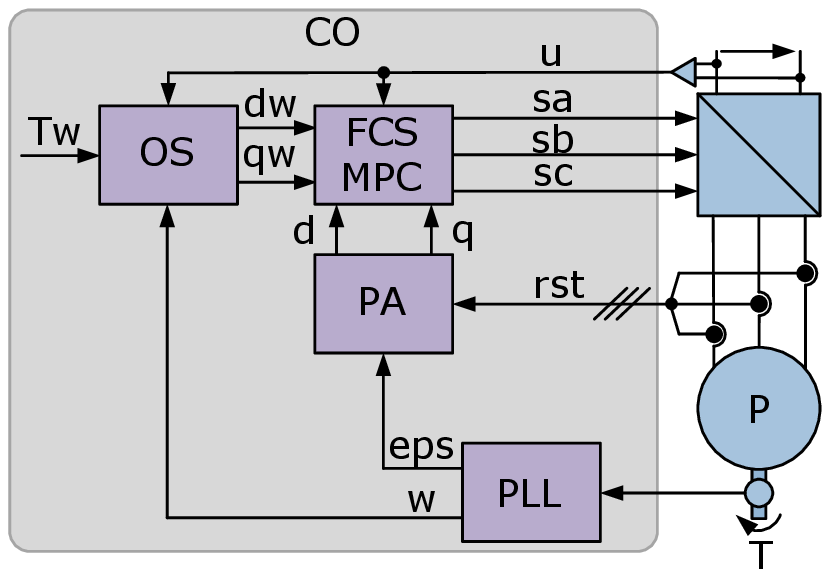}
	\else
		\includegraphics[width=0.5\columnwidth]{eps/Regelungsstruktur}
	\fi
	\caption{Field-oriented FCS-MPC topology for an IPMSM drive}
	\label{fig:Control_structure}
\end{figure}

One possible way to circumvent this trade-off decision is by using black-box plant models. This decouples the model complexity from the number of system states. 
Following the goal of streamlining the plant model, the focus of this contribution is on applying a recently developed data-driven modeling approach based on the Koopman operator \cite{PK17} to the FCS-MPC problem as depicted in \figref{fig:Control_structure}. Using this approach, dynamical systems are obtained directly for the system outputs. The interior-magnet permanent magnet synchronous motor (IPMSM) driven by a 2-level 3-phase voltage-source inverter is used as a typical power electronic application for proving the general feasibility of the Koopman operator approach. For the sake of simplicity, only the current control loop is realized by an FCS-MPC approach in this first proof of concept while the current reference values are calculated by a conventional superimposed operation point strategy (OPS). The remainder of this paper is structured as follows: In \capref{sec:MPC}, the FCS-MPC framework will be briefly summarized. The concept of the Koopman operator and its basic theoretical background will be given in \capref{sec:Koopman}. Extensive experimental tests are presented in \capref{sec:Experiment}, followed by a conclusion and an outlook in \capref{sec:Conclusion}.

\section{Finite-Set MPC Framework}
\label{sec:MPC}

Let $\bm{\Phi}$ be a discrete dynamical control system (or the flow map of a continuous system), where $\bm{x}$ is the state and $\bm{u}$ is the control input:
\begin{equation}\label{eq:ControlSystem}
\bm{x}_{i+1} = \bm{\Phi}(\bm{x}_i, \bm{u}_i), \quad i = 0,1,\ldots
\end{equation}
Here, bold symbols denote vector and matrix quantities. For the drive system given in \figref{fig:Control_structure}, the first-order Euler approximated discrete-time model in the rotor-flux oriented $dq$-system is given by:
\begin{equation}
\begin{split}
\begin{bmatrix} i_\mathrm{d,i+1}\\	i_\mathrm{q,i+1} \end{bmatrix} &= \left(\boldsymbol{I}+T_\mathrm{s}\boldsymbol{A}\right)\begin{bmatrix} i_\mathrm{d,i}\\	i_\mathrm{q,i} \end{bmatrix}+T_\mathrm{s}\bm{L}_{\mathrm{dq}}^{-1}\begin{bmatrix} u_\mathrm{d,i}\\	u_\mathrm{q,i} \end{bmatrix}+T_\mathrm{s}\begin{bmatrix} 0\\ -\frac{\psi_{\mathrm{p}}\omega_\mathrm{el}}{L_\mathrm{q}}\end{bmatrix}\\
\boldsymbol{A} &= \renewcommand*{\arraystretch}{1.4} \begin{bmatrix} -\frac{R_\mathrm{s}}{L_\mathrm{d}} & \frac{L_\mathrm{q}\omega_\mathrm{el}}{L_\mathrm{d}}\\	-\frac{L_\mathrm{d}\omega_\mathrm{el}}{L_\mathrm{q}} & -\frac{R_\mathrm{s}}{L_\mathrm{q}} \end{bmatrix} \label{eq:motor_model}, \bm{L}_{\mathrm{dq}} = \begin{bmatrix} L_\mathrm{d} & 0 \\ 0 & L_\mathrm{q}\end{bmatrix}	
\end{split}
\end{equation}
Above, $\bm{x}=\bm{i}_{\mathrm{dq}}=\begin{bmatrix} i_\mathrm{d} & i_\mathrm{q}\end{bmatrix}^\top$ is the stator current, $\bm{I}$ is the unity matrix, $T_\mathrm{s}$ is the sampling time, $\bm{L}_{\mathrm{dq}}$ is the inductance matrix, $\bm{u}_{\mathrm{dq}}=\begin{bmatrix} u_\mathrm{d} & u_\mathrm{q}\end{bmatrix}^\top$ is the stator voltage, $\psi_{\mathrm{p}}$ is the permanent magnet flux linkage, $R_\mathrm{s}$ is the stator resistance, and $\omega_\mathrm{el}$ is the electrical angular frequency, respectively.

It should be noted that temperature influences, (cross-) saturation effects and iron losses are neglected and, therefore, all motor parameters are considered constant. Linking \eqref{eq:motor_model} with \eqref{eq:ControlSystem}, the stator voltage $\bm{u}_{\mathrm{dq},i}$ is defined as:
\begin{equation*} 
	\bm{u}_{\mathrm{dq},i} = \bm{Q}(\varepsilon_{\mathrm{el},i})\frac{2}{3}\begin{bmatrix}1 & -\frac{1}{2} & -\frac{1}{2} \\[0.5em] 0 &\frac{\sqrt{3}}{2}&-\frac{\sqrt{3}}{2}\end{bmatrix}\frac{u_{\mathrm{dc},i}}{2}\bm{u}_i \
\end{equation*}
Here, $\bm{u}_i = \begin{bmatrix} s_{\mathrm{a},i} & s_{\mathrm{b},i} & s_{\mathrm{c},i}\end{bmatrix}^\top$ is the control action with the switching commands for the inverter half-bridges $s_{\mathrm{abc},i}=\left\{+1,-1\right\}$, $u_{\mathrm{dc},i}$ is the measured DC-link voltage, and $\bm{Q}$ is the rotation matrix
\begin{equation*}
\bm{Q}(\varepsilon_{\mathrm{el}}) = \begin{bmatrix} \cos(\varepsilon_{\mathrm{el}}) & \sin(\varepsilon_{\mathrm{el}}) \\[0.5em] -\sin(\varepsilon_{\mathrm{el}}) & \cos(\varepsilon_{\mathrm{el}})\end{bmatrix}	
\end{equation*}

Finally, the optimization task of the FCS-MPC for the given drive application is defined as
\begin{equation}
\begin{split}
\min_{\bm{u}_i} J = \sum_{i=1}^{n_\mathrm{p}} (i_{\mathrm{d},i}-i_{\mathrm{d},i}^*)^2 + (i_{\mathrm{q},i}-i_{\mathrm{q},i}^*)^2 
\end{split}
\label{eq:cost_fcn}
\end{equation}
with $n_\mathrm{p}$ being the number of prediction steps. To solve \eqref{eq:cost_fcn}, an exhaustive search among all possible switching sequences from $i=1, \dots, n_\mathrm{p}$ will be performed. More efficient optimization algorithms, like branch-and-bound, can improve the MPC performance because more prediction steps can be performed with the same computational effort -- however, the focus of this contribution is on the modeling part and, therefore, the exhaustive search ensures the same boundary conditions when comparing different internal FCS-MPC models.

	\section{Koopman Operator-Based Model Reduction}
\label{sec:Koopman}

In the past decade, a new approach for the construction of reduced order models (ROMs) has emerged. It is based on the Koopman operator $\Kcal$ which is a linear but infinite-dimensional operator describing the dynamics of observables of (nonlinear) autonomous dynamical systems \cite{Koo31,RMB+09}. 
By fixing the input $\bm{u}_i = \overline{\bm{u}}$ for $i = 0,1,\ldots$, the system~\eqref{eq:ControlSystem} becomes an autonomous system:
\begin{equation}\label{eq:AutonomousSystem}
\bm{x}_{i+1} = \bm{\Phi}_{\overline{\bm{u}}}(\bm{x}_i), \quad i = 0,1,\ldots
\end{equation}
where the notation $\bm{\Phi}_{\overline{\bm{u}}}$ indicates that a constant input $\overline{\bm{u}}$ is applied to system~\eqref{eq:ControlSystem}.

Introducing a real-valued observable $ \bm{f} $ of the system, 
the Koopman operator $ \Kcal_{\overline{\bm{u}}}$ (corresponding to the constant control $\overline{\bm{u}}$) describes the evolution of this observation $ \bm{y} = \bm{f}(\bm{x}) $ and is defined by
\begin{align*}
(\Kcal_{\overline{u}} \bm{f})(\bm{x}) &= \bm{f}(\bm{\Phi}_{\overline{\bm{u}}}(\bm{x}))
\end{align*}
which means that this way, a dynamical system for the observable $\bm{f}(\bm{x})$ is obtained:
Since the Koopman operator acts on observations of the dynamics, the computation is data-based, and hence, knowledge of the underlying equations is not required. Consequently, it can be used to construct a linear surrogate model from sensor data governed by nonlinear dynamics.
\begin{figure}[ht]
				\centering
				\ifnum\Review=0
					\includegraphics[width=0.8\columnwidth]{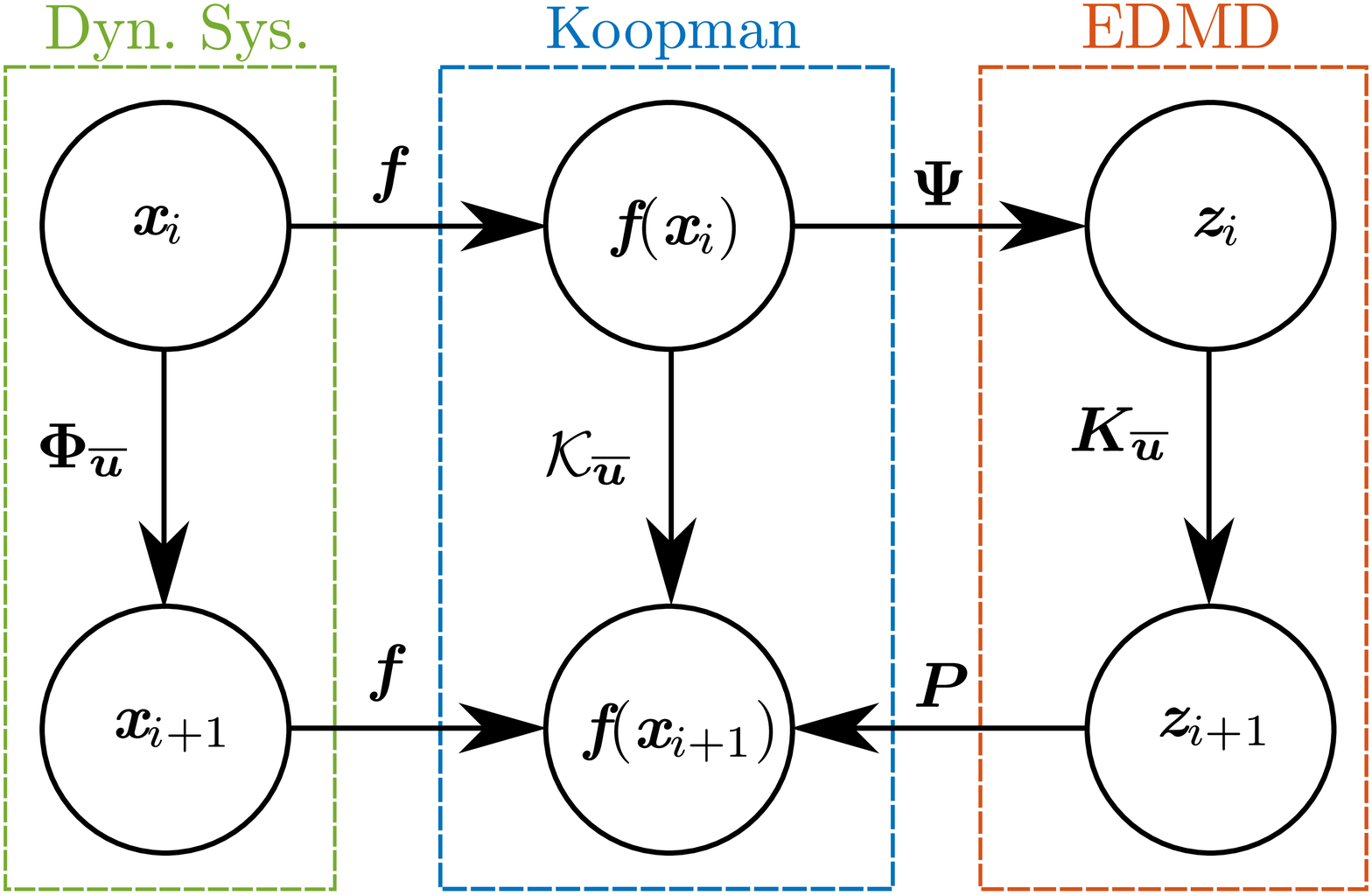}
				\else
					\includegraphics[width=0.45\columnwidth]{eps/Koopman_Diagramm5}
				\fi
				\caption{Relation between the system dynamics $\Phi_{\overline{u}}$, the corresponding Koopman operator $\Kcal_{\overline{u}}$ and its finite-dimensional representation $K_{\overline{u}}$ computed via EDMD}
				\label{fig:Koopman_EDMD}
			\end{figure}
			
The most popular approach to construct a finite-dimensional approximation of the Koopman operator is via Dynamic Mode Decomposition (DMD) \cite{Sch10} or Extended Dynamic Mode Decomposition (EDMD) \cite{WKR15}. 
While in DMD the Koopman operator is approximated for the observations themselves, in EDMD they are expressed in terms of arbitrary basis functions (e.g., monomials, Hermite polynomials or radial basis functions). 
For a given set of basis functions $ \{ {\psi}_{1},\,{\psi}_{2},\,\dots,\,{\psi}_{k} \} $ (the so-called dictionary), define $ \bm{\psi} $ by
\begin{equation*}
\bm{\psi}(\bm{y}) =
\left[{\psi}_{1}(\bm{y}) \quad {\psi}_{2}(\bm{y}) \quad \dots \quad {\psi}_{k}(\bm{y})\right]^{\top},
\end{equation*}
where  $\bm{y} = \bm{f}(\bm{x})$. Choosing $ \bm{\psi}(\bm{y}) = \bm{y} $ yields DMD as a special case of EDMD. For the computation, either measurement or simulation data is used which is written in matrix form as
\begin{equation*}
\bm{Y} =
\begin{bmatrix}
\bm{y}_1 & \bm{y}_2 & \cdots & \bm{y}_m
\end{bmatrix}
\quad\text{and}\quad
\widehat{\bm{Y}} =
\begin{bmatrix}
\widehat{\bm{y}}_1 & \widehat{\bm{y}}_2 & \cdots & \widehat{\bm{y}}_m
\end{bmatrix}
\end{equation*}
where $ \widehat{\bm{y}}_i = \bm{f}(\bm{\Phi}_{\overline{\bm{u}}}(\bm{x}_i)) $. The data can be obtained from one long trajectory such that $ \widehat{\bm{y}}_i = \bm{y}_{i+1} $. Alternatively, many short simulations may be assembled in the matrices $\bm{Y}$ and $\widehat{\bm{Y}}$. For EDMD, the data matrices are embedded into the typically higher-dimensional feature space by
\begin{align*}
\bm{\Psi}_{\bm{Y}} &=\left[ \bm{\psi}(\bm{y}_{1}) \quad \bm{\psi}(\bm{y}_{2}) \quad \dots \quad \bm{\psi}(\bm{y}_{m}) \right]\\
\bm{\Psi}_{\widehat{\bm{Y}}} &=
\left[ \bm{\psi}(\widehat{\bm{y}}_{1}) \quad \bm{\psi}(\widehat{\bm{y}}_{2}) \quad \dots \quad \bm{\psi}(\widehat{\bm{y}}_{m}) \right]
\end{align*}

With these data matrices, the matrix $ \bm{K}_{\overline{\bm{u}}} \in \R^{k \times k} $ can be computed, where $k$ depends on the dimension of the dictionary:
\begin{equation}
\bm{K}_{\overline{\bm{u}}}^{\top} = \bm{\Psi}_{\widehat{Y}} \bm{\Psi}_Y^+ = \big( \bm{\Psi}_{\widehat{Y}} \bm{\Psi}_Y^{\top} \big) \big(\bm{\Psi}_Y \bm{\Psi}_Y^{\top}\big)^+,
\label{eq:koopman_aprox}
\end{equation}
see \cite{KKS16} for details. The matrix $ \bm{K}_{\overline{\bm{u}}} $ can be viewed as a finite-dimensional approximation of the Koopman operator. 
Instead of the more common approach, where a decomposition into eigenvalues, eigenfunctions, and modes is applied to analyze the system dynamics, 
updates for the observable $\bm{f}(\bm{x})$ can be computed directly using $\bm{K}_{\overline{\bm{u}}}$, cf.~\cite{PK17}:
\begin{equation}
\bm{z}_{i+1} = \bm{K}_{\overline{\bm{u}}}^{\top} \bm{z}_{i}, \quad i = 0,1,\ldots, \label{eq:K-ROM}
\end{equation}
where $\bm{z} = \bm{\psi}(\bm{f}(\bm{x}))$.
This means that a linear surrogate model for a potentially nonlinear and unknown dynamical system is constructed from data.
From here, $\bm{f}(\bm{x}_{i+1})$ can be obtained using the projection matrix $\bm{P}$, cf.~\figref{fig:Koopman_EDMD}, where the relation between the dynamical system $\bm{\Phi}_{\overline{\bm{u}}}$, the related Koopman operator $\Kcal_{\overline{\bm{u}}}$, and the EDMD approximation $\bm{K}_{\overline{\bm{u}}}$ is visualized. Note that $\bm{z} = \bm{f}(\bm{x})$ if DMD is used instead of EDMD, which will be done in this article.

More recently, various attempts have been made to use ROMs based on the Koopman operator for both open and closed-loop control problems \cite{PBK15,KM18a,KKB17}. In these approaches, the Koopman operator is either approximated for an augmented state (consisting of the actual state and the control) in order to deal with the non-autonomous control system or an affine control dependency is assumed. An alternative approach which will be utilized in this article is to replace the control system \eqref{eq:ControlSystem} by a set of autonomous systems \eqref{eq:AutonomousSystem} with constant control inputs $\overline{\bm{u}}_1, \ldots, \overline{\bm{u}}_{n}$. This way, the optimal control problem is either transformed into a switching time problem \cite{PK17} or into a bilinear control problem \cite{Pei18}. 

Since the finite-control-set MPC problem introduced in \capref{sec:MPC} is already a problem of switching type, the approach from \cite{PK17} can be applied directly. The only adaptation necessary is to replace the original dynamics \eqref{eq:motor_model} by Koopman operator-based ROMs \eqref{eq:K-ROM}. In the example given in \figref{fig:Control_structure}, the 2-level inverter has $2^3=8$ switching states, but the two zero-voltage vectors lead to the same autonomous system behavior. Even though inverter power losses may vary depending on which of the two zero-voltage inputs is applied, this will be neglected for the sake of simplicity. Hence, only seven Koopman operator-based ROMs have to be computed which reduces the online computational burden for the FCS-MPC due to fewer switching sequences.
\begin{remark}
	For the problem considered here, a linearization can be obtained very efficiently, cf.~\eqref{eq:motor_model}. In this situation, the Koopman operator-based approach is similar to a simple least squares fit when using DMD. However, it should be noted that the method presented here does not require a linear or linearized system but can be used to construct a linear model for highly nonlinear and even unknown system dynamics.
\end{remark}

In an offline phase, data is collected for all seven system states and the matrices $\bm{K}_{\overline{\bm{u}}_1}$ to $\bm{K}_{\overline{\bm{u}}_7}$ are constructed. Since real-time applicability is crucial, DMD is used instead of EDMD such that the matrices have a very low dimension. Furthermore, numerical experiments have shown that this is the most robust approach in the low data limit \cite{PK18}.
The relevant observations are the currents $i_\mathrm{d}$ and $i_\mathrm{q}$ as well as the electrical rotor angle $\varepsilon_\mathrm{el}$. Nearly constant rotational speed over the prediction horizon is assumed, i.e. the rotor angle is increasing continuously such that the resulting Koopman operator is unbounded which violates the requirements for the convergence properties and yields unsatisfactory results \cite{KM18}. To this end, the sine and cosine of the angle are observed which results in the observation \eqref{eq:observation}. The corresponding ROM hence has dimension four.
\begin{equation}
	\bm{y} = \left[ i_\mathrm{d} \quad i_\mathrm{q} \quad \sin(\varepsilon_{\mathrm{el}}) \quad \cos(\varepsilon_{\mathrm{el}}) \right]^{\top}
\label{eq:observation}
\end{equation}

\begin{figure}[ht]
	\centering
	\psfrag{Mdl}[cc][cc][1][0]{Linear motor model \eqref{eq:motor_model}}
	\psfrag{Kop}[cc][cc][1][0]{Koopman matrix training}
	\psfrag{Si}[cc][cc][0.95][0]{Offl. simulation}
	\psfrag{Opt}[cc][cc][0.95][0]{ROM fitting \eqref{eq:koopman_aprox}}
	\psfrag{FCS}[cc][cc][1][0]{FCS-MPC}
	\psfrag{Y}[cc][cc][1][0]{$\hat{\bm{Y}}$}
	\psfrag{y}[cc][cc][1][0]{$\bm{Y}$}
	\psfrag{eq2}[cc][cc][1][0]{$\bm{\Phi}$}
	\psfrag{eq1}[cc][cc][1][0]{$\bm{K}_{\overline{\bm{u}}}^{\top}$}
		\ifnum\Review=0
		\includegraphics[width=0.99\columnwidth]{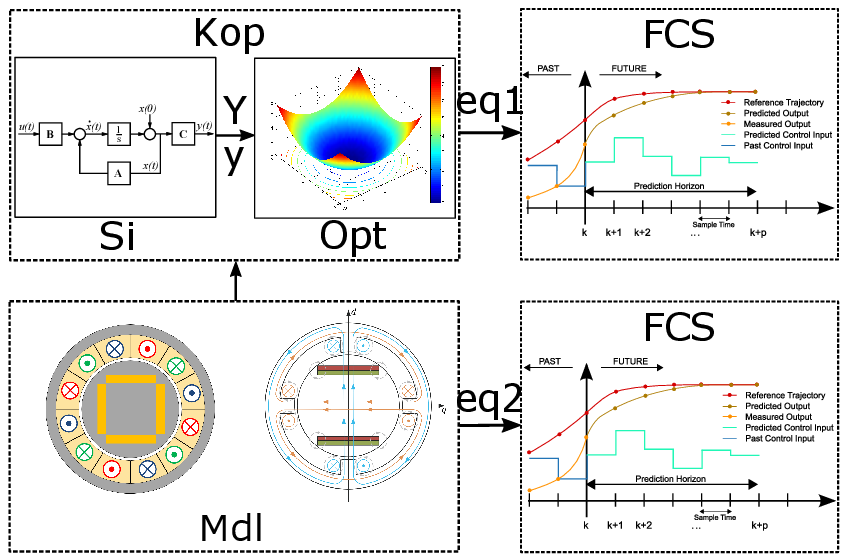}
	\else
		\includegraphics[width=0.65\columnwidth]{eps/Koopman_Training}
	\fi
	\caption{Koopman operator-based ROM training and FCS-MPC application}
	\label{fig:Koopman_Training}
\end{figure}

For the following experimental validation, the observations have been generated according to \figref{fig:Koopman_Training}: The time-continuous variant of the standard white-box motor model \eqref{eq:motor_model} has been utilized in a simple Simulink-based closed-loop control simulation to generate data linking the control action $\bm{u}_i$ to the observation $\bm{y}$. Using this data, the matrices $\bm{K}_{\overline{\bm{u}}_1}$ to $\bm{K}_{\overline{\bm{u}}_7}$ are computed via \eqref{eq:koopman_aprox}. 

Following this simulation-based offline training process, it becomes clear that the computed Koopman approximation can perform, at best, at the same level compared to the original model \eqref{eq:motor_model} in an FCS-MPC application. The use of measurement (i.e., model-free) training data is presently being investigated. In addition, it should be pointed out that the training was only carried out for one fixed motor speed of \unit{\textit{n}=1000}{min^{-1}}. Operating the Koopman-MPC at other speeds will result in a systematic modeling error.  This issue will be discussed in \capref{subsec:param_distort}. Thus, the primary objectives of this publication are to evaluate the approximation accuracy of the Koopman-based ROM as well as its computational load compared to the baseline model and also its general feasibility in a MPC context. 

\section{Experimental Evaluation}
\label{sec:Experiment}

The proposed Koopman operator-based FCS-MPC has been implemented on a laboratory test bench equipped with a dSPACE DS1006 rapid-control-prototyping system. The drive system under test consists of a \unit{55}{kW} IPMSM manufactured by Brusa (HSM1-6.17.12-C01) and a 2-level IGBT inverter from Semikron (SKiiP 1242GB120-4D). Moreover, a speed-controlled load machine is coupled with the DUT via torque meter. The most important motor and test bench parameters are summarized in \tabref{tab:parm}.
	\begin{table}[ht]
				\renewcommand{\arraystretch}{1.15}
				\caption{IPMSM and control Parameters}
			 \label{tab:parm}
				\centering
				\begin{tabular}{l|c|c}
				\hline
				Stator resistance         & $R$               & 18\,m$\Omega$\\
				Inductance in d-direction & $L_\mathrm{d}$    & 370\,$\micro$H\\
				Inductance in q-direction & $L_\mathrm{q}$    & 1200\,$\micro$H\\
				Permanent magnet flux     & $\Psi_\mathrm{pm}$ & 66\,mVs\\ 
				Pole pair number          & $p$               & 3\\ 
				\hline
				DC-link voltage                  		& $u_\mathrm{DC}$   							  	& 300\,V\\
				Mechanical speed                 		& $n$   													  	& 1000\,min$^{-1}$\\
				FOC: \,Controller cycle time    	  & $T_\mathrm{s,FOC}$								  & 50\,$\micro$s \\
				FOC: \,Symmetrical optimum 					& $a$               								  & 3 or 4\\
        FOC: \,Oversampling factor      		& $f_\mathrm{s,FOC}/f_\mathrm{sw}$    & 5 or 6\\
				MPC: Controller cycle time       		& $T_\mathrm{s,MPC}$  					  		& 50$\,\micro$s\\
				MPC: Max. switching frequency    		& $f_\mathrm{sw}$   						  		& 10\,kHz\\
				MPC: Prediction horizon    		      & $n_\mathrm{p}$   	  					  		& 3\\
				\hline
				\end{tabular}
		\end{table}	
					
At the test bench, three different control approaches has been implemented and tested:
\begin{itemize}
	\item FCS-MPC with Koopman operator-based ROM
	\item FCS-MPC with standard white-box motor model \eqref{eq:motor_model}
	\item Field-oriented control (FOC) with PWM
\end{itemize}
Both MPC variants compensate for the one step delay due to the digital implementation of the controllers by adding a prediction step before starting the exhaustive search with $n_p=3$ steps. The FCS-MPC approach in conjunction with the selected controller cycle time of 50$\,\micro$s results in an upper limit of the switching frequency of 10\,kHz. The FOC was analytically tuned according to the well-known symmetrical optimum. For a better comparability between FOC and the MPC variants, the FOC was implemented in an oversampling framework \cite{BoeckerBuchholz2013}. In this way, the controller operates at the same cycle time compared to the MPCs while the frequency of the carrier is manually adapted to the resulting switching frequencies of the MPC realizations. 

The turnaround times of the implemented controllers on the used hardware setup are summarized in \tabref{tab:turnaround_times}. Among other test bench specific computational overhead, the time needed for the analog to digital conversions and the PLL are included. It can be seen that the Koopman-based execution times are in the same range compared to the white-box model-based MPC. It can be concluded that the simplest white-box motor model has the same computational complexity as the Koopman-based model. When additional effects, like (cross-)  saturation or iron losses, should be modeled within the white-box approach, this will directly result in increasing computational load, whereas the Koopman-based black-box model can be trained with the same Koopman matrix dimensions and, therefore, preserve the computational complexity regarding online computations. 
\begin{table}[ht]
										\renewcommand{\arraystretch}{1.15}
										\caption{Turnaround times of the control methods}
										\label{tab:turnaround_times}
										\centering
										\begin{tabular}{l|c|c}                   
																 & mean value & standard deviation\\ \hline
										FOC          & \unit{27.6}{\mu s}		& \unit{0.42}{\mu s}  \\
										Koopman-MPC  & \unit{28.5}{\mu s}  & \unit{0.40}{\mu s}   \\
										Standard-MPC   & \unit{26.5}{\mu s}  & \unit{0.39}{\mu s}    \\ 
										\hline
										\end{tabular}
		\end{table}

\subsection{Evaluation at nominal speed}

\begin{figure}[ht]
			\centering
				\includegraphics{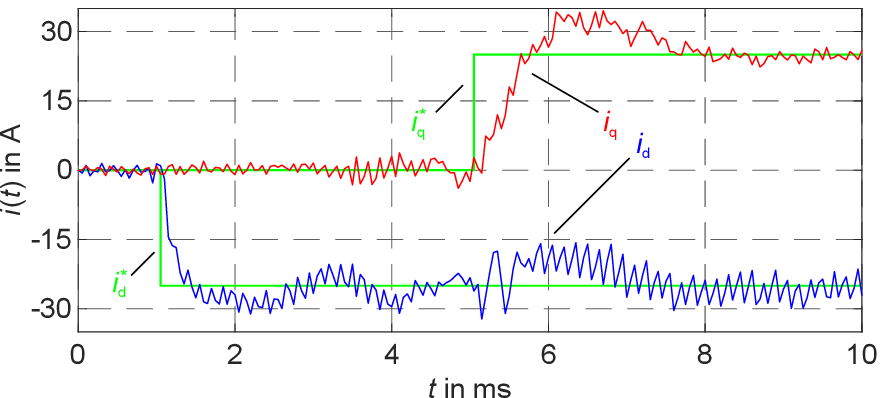}
				\caption{FOC: small signal response at $n$=\unit{1000}{min^{-1}}, overs. factor 6, $a$=3}
				\label{FOC_small_signal}
				\end{figure}

\begin{figure}[h!]
			\centering
				\includegraphics{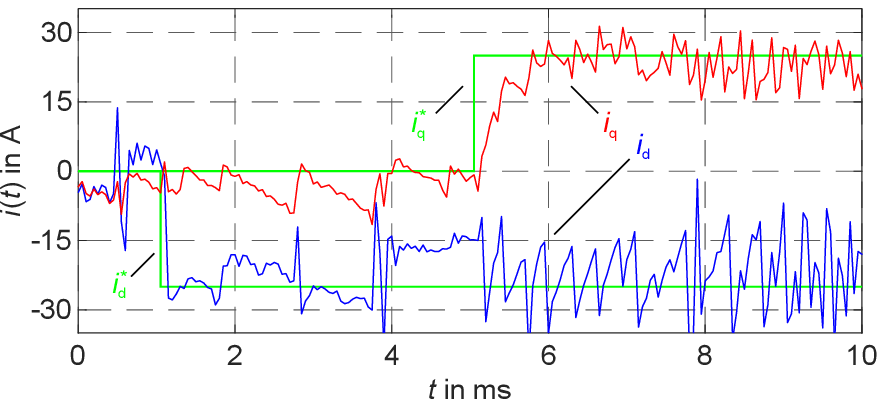}
				\caption{Koopman-MPC: small signal response at $n$=\unit{1000}{min^{-1}}}
				\label{Koopman_small_signal}
				\end{figure}

\begin{figure}[h!]
			\centering
				\includegraphics{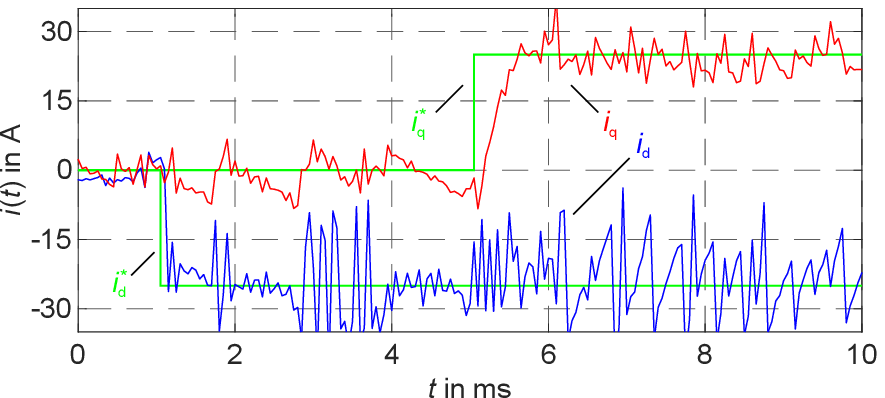}
				\caption{Standard-MPC: small signal response at $n$=\unit{1000}{min^{-1}}}
				\label{Linear_small_signal}
				\end{figure}

For the comparison of the three considered control approaches, the step response as well as the behavior in steady state at the nominal fixed rotational speed of \unit{\textit{n}=1000}{min^{-1}} were examined. The responses to small current reference steps are shown in \figref{FOC_small_signal}, \figref{Koopman_small_signal} and \figref{Linear_small_signal}, respectively. The reference $i^{*}_\mathrm{d}$ is changed to -25\,A first. Afterwards $i^{*}_\mathrm{q}$ is set to 25\,A. For the shown responses, the currents sampled by the respective controller are used. The carrier frequency of the PWM used for the FOC is set to 3.3\,kHz, resulting in an oversampling factor between PI-controller and carrier of six. The controller is tuned with $a=3$ for the symmetrical optimum.

\setlength\abovecaptionskip{-10pt}
\setlength\abovecaptionskip{6pt}
	
Both MPC variants are able to comply with the reference change. The behavior is very similar. Compared to the FOC, the settling times are significantly shorter. While the FOC needs approximately 3\,ms after the reference step for the q-current, the MPCs need less than 1\,ms. Also, the decoupling between the two current axes is slightly better when using MPC compared to FOC. For each control approach, the total harmonic distortion (THD) of the phase currents and the setpoint deviations in steady state operation are summarized in \tabref{tab:THD_setpoint_deviation}. The setpoint deviation is evaluated as the geometric deviation between a sliding mean of the d- and q-current and the corresponding setpoint. The average switching frequencies for different steady state operations are summarized in \tabref{tab:switching_frequencies}.
 
All three approaches use an average switching frequency of 3.3\,kHz. The current ripple and the THD with the FOC are lower than with the MPCs, which directly follows from the degree of freedom to quasi-continuously shift the switchovers in time due to the used PWM. For the FCS-MPCs, switchovers can take place only at discrete time points. Due to the used PI-controller, the setpoint deviation using the FOC is nearly zero. In contrast, the deviation with the MPCs is approximately 2.5\,A.

\begin{table}[ht]
				\renewcommand{\arraystretch}{1.15}
				\caption{THD and setpoint deviation at $n$=\unit{1000}{min^{-1}}}
				\label{tab:THD_setpoint_deviation}
				\centering
				\begin{tabular}{l|c|c|c|c}                   
										 & \multicolumn{2}{c|}{$i_\mathrm{d}$ = -25\,A, $i_\mathrm{q}$ = 25\,A}		& \multicolumn{2}{c}{$i_\mathrm{d}$ = -169\,A, $i_\mathrm{q}$ = 169\,A} \\ 
										 & THD   & setp. deviation                & THD      & setp. deviation  \\ \hline
				FOC          & 8.7\,\% & 0.0\,A                         & 4.7\,\%  & 0.0\,A      \\
				Koopman-MPC  & 21.1\,\% & 2.6\,A                           &	15.6\,\% & 1.2\,A    \\
				Standard-MPC & 20.5\,\% & 2.3\,A   		                  &	15.2\,\% & 1.2\,A      \\
				\hline
				\end{tabular}
		\end{table} 
\begin{table}[ht]
				\renewcommand{\arraystretch}{1.15}
				\caption{Average switching frequencies at $n$=\unit{1000}{min^{-1}}}
				\label{tab:switching_frequencies}
				\centering
				\begin{tabular}{l|c|c}                   
										 & $i_\mathrm{d}$ = -25\,A, $i_\mathrm{q}$ = 25\,A & $i_\mathrm{d}$ = -169\,A, $i_\mathrm{q}$ = 169\,A \\ \hline
				FOC          & 3.3\,kHz																		     & 4.0\,kHz \\
				Koopman-MPC  & 3.3\,kHz                                        & 4.0\,kHz \\
				Standard-MPC & 3.3\,kHz                                        & 4.0\,kHz \\
				\hline
				\end{tabular}
		\end{table}

\begin{figure}[h!]
			\centering
				\includegraphics{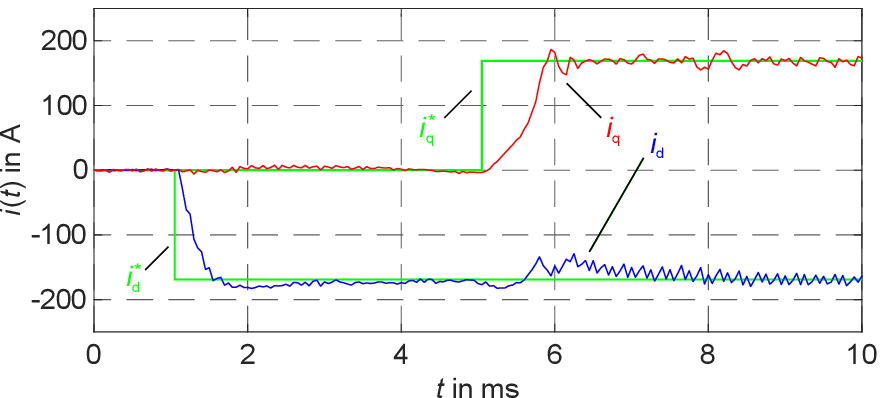}
				\caption{FOC: large signal response at $n$=\unit{1000}{min^{-1}}, overs. factor 5, $a$=4} 
				\label{FOC_large_signal}
		\end{figure}

\begin{figure}[h!]
			\centering
				\includegraphics{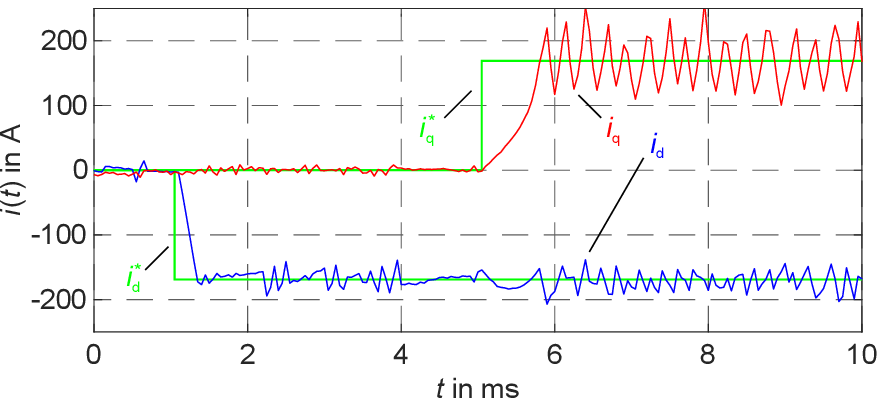}
				\caption{Koopman-MPC: large signal response at $n$=\unit{1000}{min^{-1}}} 
				\label{Koopman_large_signal}
		\end{figure}

\begin{figure}[h!]
			\centering
				\includegraphics{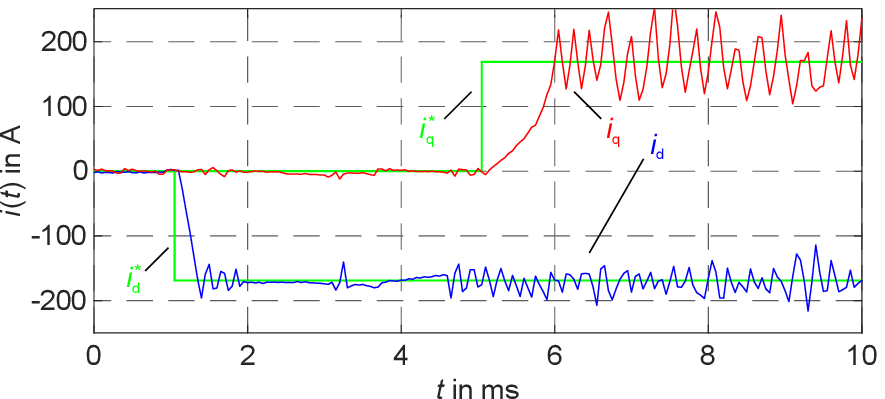}
				\caption{Standard-MPC: large signal response at $n$=\unit{1000}{min^{-1}}} 
				\label{Linear_large_signal}
		\end{figure}

\figref{FOC_large_signal}, \figref{Koopman_large_signal} and \figref{Linear_large_signal} show the responses to large current reference steps. Here, the currents are changed from zero to $i_\mathrm{d}^*$ = -169\,A and $i_\mathrm{q}^*$ = 169\,A which corresponds to the nominal current of the motor. For the FOC, an oversampling factor of five and $a=4$ are chosen. After the setpoint change, the d-current reaches the steady state operation at about 1.5\,ms. In less than 1\,ms the set point is reached with both MPCs. The rise time of the q-current is nearly identical in all cases. During the step responses of the MPCs there is hardly any mutual influencing between the two currents. Using the FOC, a significant influence on the d-current during the rise of the q-current can be noticed. Due to saturation effects the ripple increases at higher currents.

\begin{figure}[ht]
			\centering
				\includegraphics{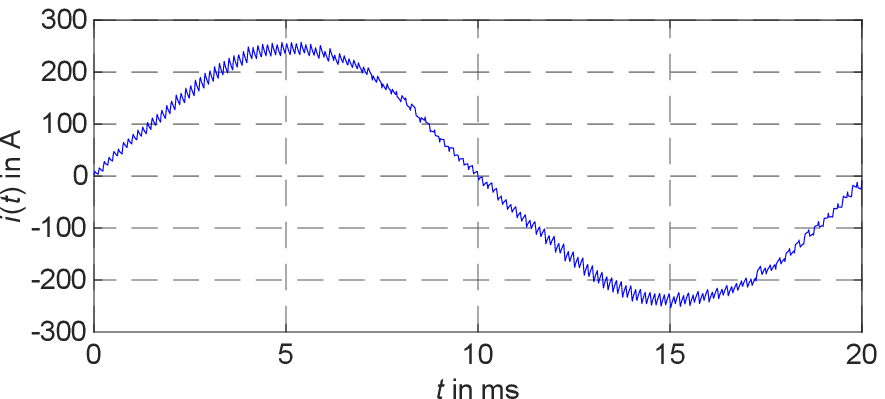}
				\caption{FOC: $i_\mathrm{a}$ at $n$=\unit{1000}{min^{-1}}, overs. factor 5, $a$=4}
				\label{FOC_phase_current}
		\end{figure}

\begin{figure}[h!]
			\centering
				\includegraphics{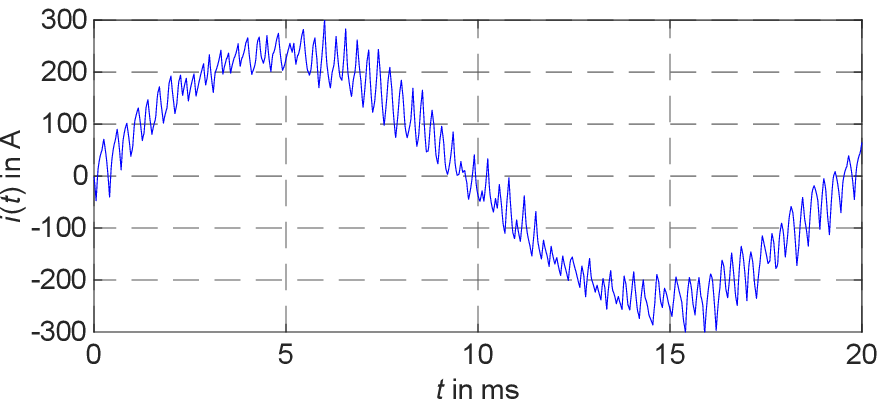}
				\caption{Koopman-MPC: $i_\mathrm{a}$ at $n$=\unit{1000}{min^{-1}}}
				\label{Koopman_phase_current}
		\end{figure}

\begin{figure}[h!]
			\centering
				\includegraphics{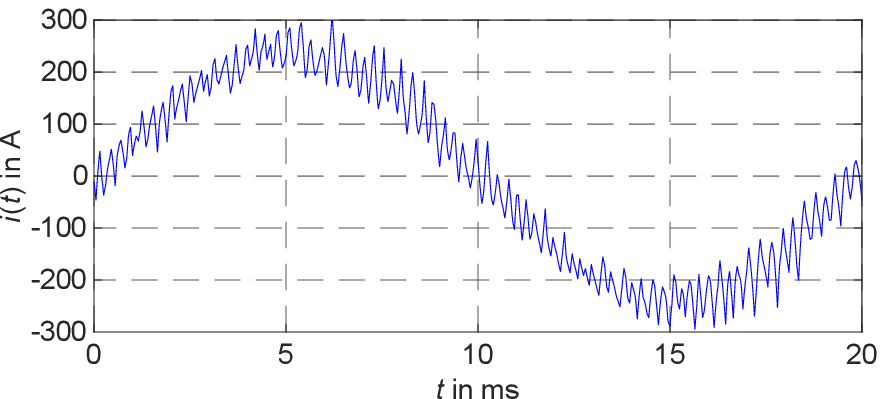}
				\caption{Standard-MPC: $i_\mathrm{a}$ at $n$=\unit{1000}{min^{-1}}}
				\label{Linear_phase_current}
		\end{figure}

\begin{figure}[ht]
			\centering
				\includegraphics{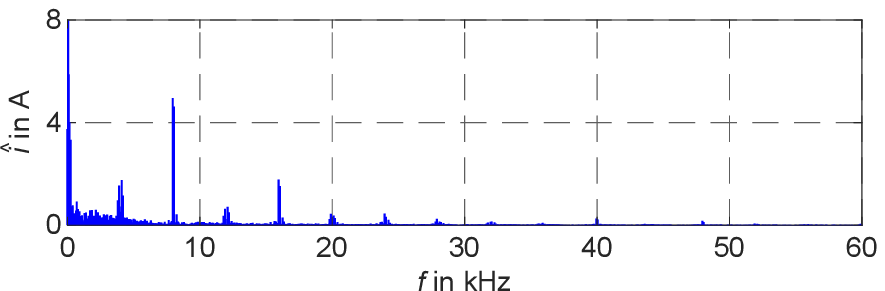}   
				\caption{FOC: DFT of $i_\mathrm{a}$ at $n$=\unit{1000}{min^{-1}}, overs. factor 5, $a$=4}
				\label{FOC_DFT}
		\end{figure}

\begin{figure}[h!]
			\centering
				\includegraphics{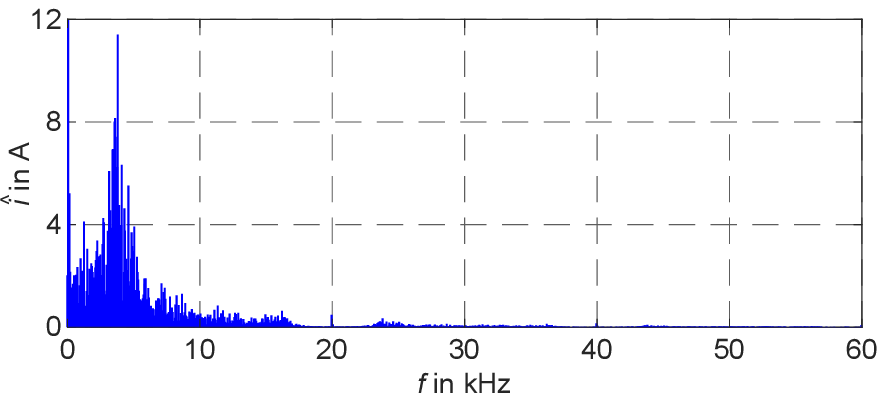}
				\caption{Koopman-MPC: DFT of $i_\mathrm{a}$ at $n$=\unit{1000}{min^{-1}}}
				\label{Koopman_DFT}
		\end{figure}

\begin{figure}[h!]
			\centering
				\includegraphics{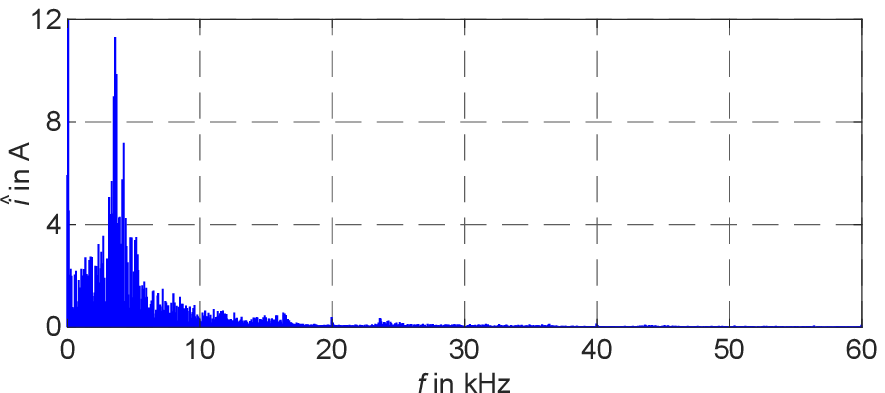}
				\caption{Standard-MPC: DFT of $i_\mathrm{a}$ at $n$=\unit{1000}{min^{-1}}}
				\label{Linear_DFT}
		\end{figure}

The currents of one phase during operation with the nominal current operating points are shown in  \figref{FOC_phase_current}, \figref{Koopman_phase_current} and \figref{Linear_phase_current} to determine the current ripple and harmonics. The currents are sampled at a rate of 50\,MHz by means of an external transient recorder. The corresponding discrete Fourier transformations (DFT) are given in \figref{FOC_DFT}, \figref{Koopman_DFT} and \figref{Linear_DFT}. The average switching frequency for all three controllers is 4.0\,kHz. 

For the FOC, current harmonics at multiples of the switching frequency emerge in the spectra, and the resulting total harmonic distortion (THD) is around 5\,\%, whereas the THD with the MPCs is around 15\,\%. Here, the higher-order harmonics are located mainly in a broader spectrum around the mean switching frequency. The setpoint deviation with the MPCs is approximately 1\,A.

\setlength\abovecaptionskip{-10pt}
\setlength\abovecaptionskip{6pt}

\setlength\abovecaptionskip{-10pt}
\setlength\abovecaptionskip{6pt}

\subsection{Evaluation at deviating speeds}
\label{subsec:param_distort}

The Koopman-based ROM was constructed from data obtained at a constant mechanical speed of \unit{\textit{n}=1000}{min^{-1}}. To investigate the robustness of the ROM performance, the effects of deviating speeds are considered. For comparison, the white-box standard model-based FCS-MPC is utilized. The latter involves updating the variables in the model, such as the speed, at the beginning of the next controller cycle. Since the typical behavioral differences compared to the FOC approach are similar to the previous section, it will not be discussed in this section again. 

\begin{figure}[ht]
			\centering
				\includegraphics{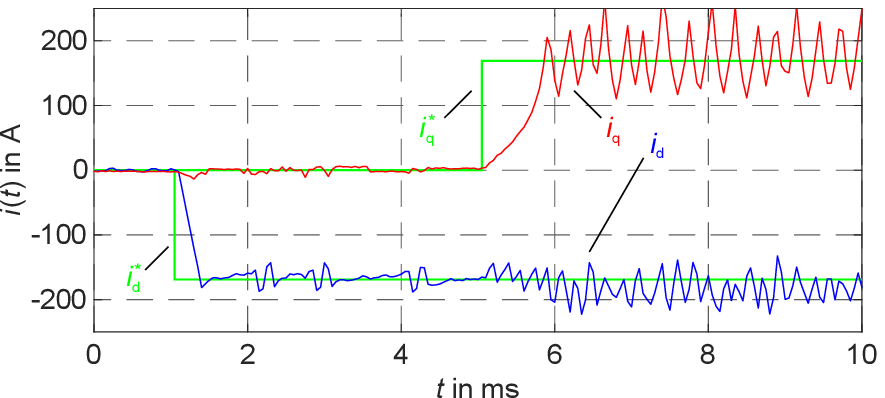}
				\caption{Koopman-MPC: large signal response at $n$=\unit{100}{min^{-1}}} 
				\label{Koopman_large_signal_n100}
		\end{figure}

\begin{figure}[h!]
			\centering
				\includegraphics{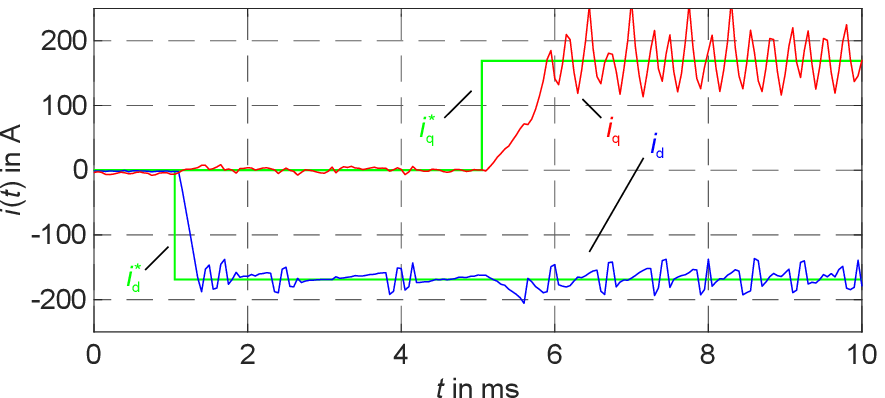}
				\caption{Standard-MPC: large signal response at $n$=\unit{100}{min^{-1}}} 
				\label{Linear_large_signal_n100}
		\end{figure}

\begin{figure}[h!]
			\centering
				\includegraphics{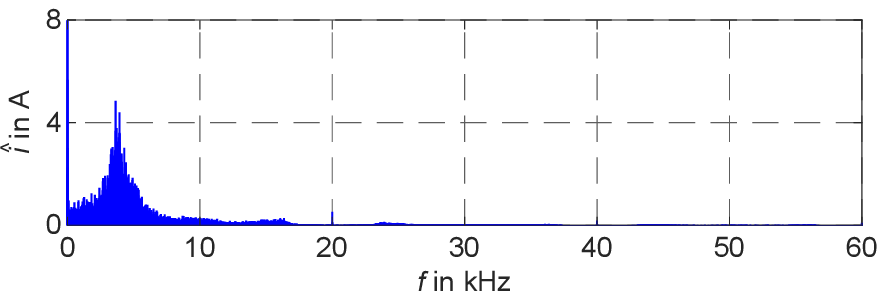}    
				\caption{Koopman-MPC: DFT of $i_\mathrm{a}$ at $n$=\unit{100}{min^{-1}} and nominal current}
				\label{Koopman_DFT_n100}
		\end{figure}

\begin{figure}[h!]
			\centering
				\includegraphics{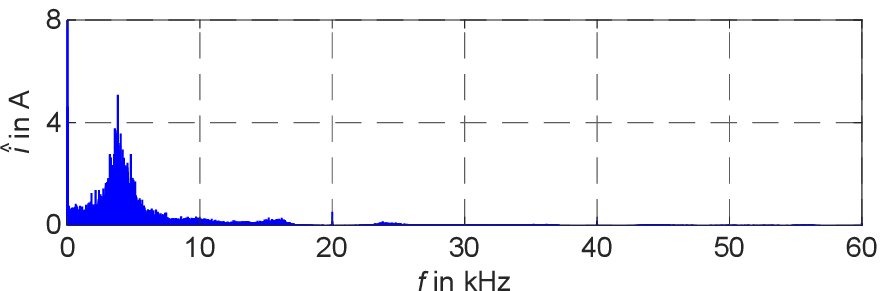}   
				\caption{Standard-MPC: DFT of $i_\mathrm{a}$ at $n$=\unit{100}{min^{-1}} and nominal current}
				\label{Linear_DFT_n100}
		\end{figure}

In \figref{Koopman_large_signal_n100} and \figref{Linear_large_signal_n100} the large signal step responses for both MPC approaches are shown for a mechanical speed of \unit{\textit{n}=100}{min^{-1}}. Moreover, the corresponding current spectra during steady state at nominal current are shown in \figref{Koopman_DFT_n100} and \figref{Linear_DFT_n100}. It can be seen that both modeling approaches deliver nearly equal dynamic behavior in the FCS-MPC framework. However, the setpoint deviation of the Koopman-MPC increases considerably compared to the Standard-MPC. The specific deviations for both speeds are given in \tabref{tab:setpoint_deviation_deviating_speeds}, the average switching frequencies in \tabref{tab:switching_frequencies_deviation_deviating_speeds}. 

\begin{table}[ht]
				\renewcommand{\arraystretch}{1.15}
				\caption{THD and setpoint deviation at deviating speeds,\newline $i_\mathrm{d}$=-169\,A, $i_\mathrm{q}$=169\,A}
				\label{tab:setpoint_deviation_deviating_speeds}
				\centering
				\begin{tabular}{l|c|c|c|c}                   
										 & \multicolumn{2}{c|}{\unit{100}{min^{-1}}}		& \multicolumn{2}{c}{\unit{2500}{min^{-1}}} \\ 
										 & THD      & setp. deviation                   & THD      & setp. deviation  \\ \hline
				Koopman-MPC  & 14.2\,\% & 14.3\,A                           &	14.3\,\% & 20.3\,A      \\
				Standard-MPC & 15.5\,\% & 3.0\,A      		                  &	13.5\,\% & 5.2\,A       \\
				\hline
				\end{tabular}
		\end{table} 
\begin{table}[ht]
				\renewcommand{\arraystretch}{1.15}
				\caption{Average switching frequencies at deviating speeds, $i_\mathrm{d}$=-169\,A, $i_\mathrm{q}$=169\,A}
				\label{tab:switching_frequencies_deviation_deviating_speeds}
				\centering
				\begin{tabular}{l|c|c}                   
										 & \unit{100}{min^{-1}}		& \unit{2500}{min^{-1}} \\ \hline
				Koopman-MPC  & 4.2\,kHz   						& 2.8\,kHz  \\
				Standard-MPC & 4.3\,kHz    						& 2.7\,kHz  \\
				\hline
				\end{tabular}
		\end{table}
		
\setlength\abovecaptionskip{-10pt}		
\setlength\abovecaptionskip{6pt}

\begin{figure}[ht]
			\centering
				\includegraphics{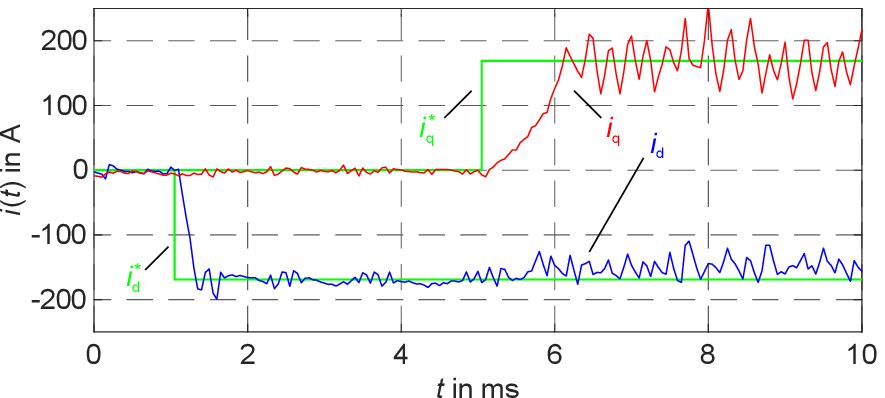}
				\caption{Koopman-MPC: large signal response at $n$=\unit{2500}{min^{-1}}} 
				\label{Koopman_large_signal_n2500}
		\end{figure}

\begin{figure}[h!]
			\centering
				\includegraphics{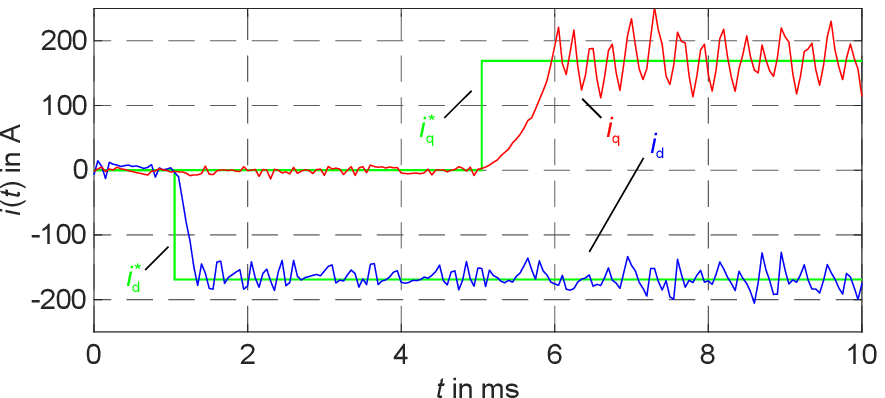}
				\caption{Standard-MPC: large signal response at $n$=\unit{2500}{min^{-1}}} 
				\label{Linear_large_signal_n2500}
		\end{figure}

\begin{figure}[h!]
			\centering
				\includegraphics{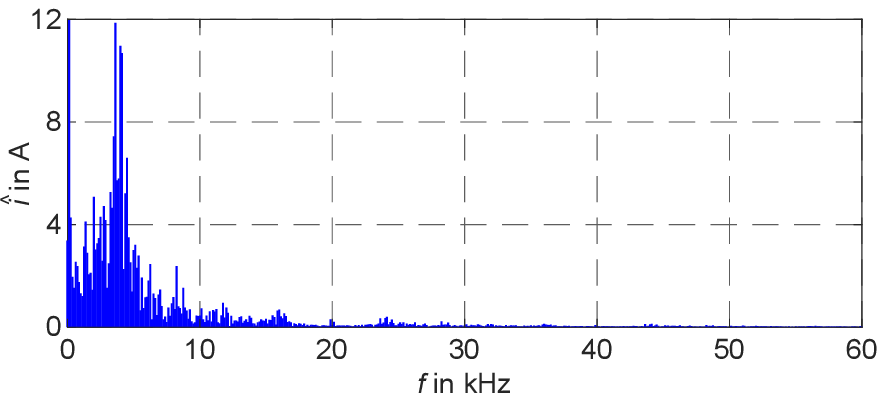}
				\caption{Koopman-MPC: DFT of $i_\mathrm{a}$ at $n$=\unit{2500}{min^{-1}} and nominal current}
				\label{Koopman_DFT_n2500}
		\end{figure}

\begin{figure}[h!]
			\centering
				\includegraphics{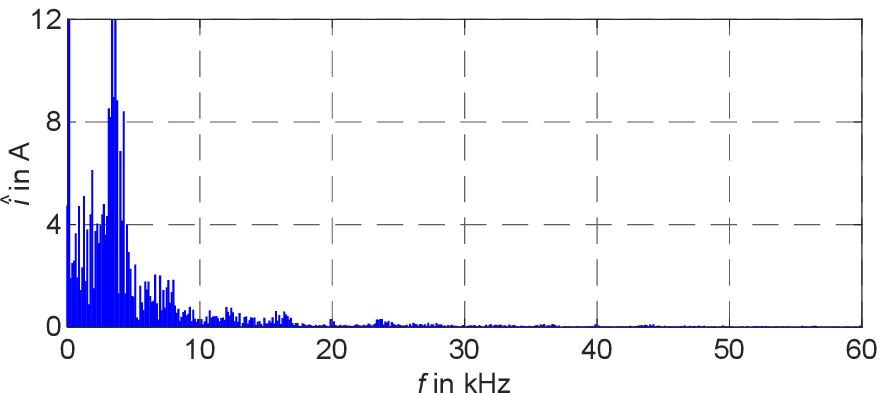}
				\caption{Standard-MPC: DFT of $i_\mathrm{a}$ at $n$=\unit{2500}{min^{-1}} and nominal current}
				\label{Linear_DFT_n2500}
		\end{figure}

Furthermore, the control behavior is also evaluated at a higher speed of \unit{\textit{n}=2500}{min^{-1}}. The large signal step responses are depicted in \figref{Koopman_large_signal_n2500} and \figref{Linear_large_signal_n2500} while the current spectra at steady state are given in \figref{Koopman_DFT_n2500} and \figref{Linear_DFT_n2500}. Again, the dynamic behavior of both plant modeling approaches is almost identical. Here as well, the setpoint deviation for the Koopman-MPC shows a strong increase that results from the chosen structure of the ROM, which does not take into account the influences of a deviating speed. However, through an online update based on measurement data, these influences can be taken into consideration.

\setlength\abovecaptionskip{-10pt}
\setlength\abovecaptionskip{6pt}

\section{Conclusion and Outlook}
\label{sec:Conclusion}
Extensive experimental tests have proven that the Koopman operator-based FCS-MPC achieves the same performance as the standard white-box model-based approach which highlights the robustness of the presented black-box plant model. Even when tested at significantly different operation regimes compared to the training setup, the presented method shows satisfying results. To the authors' best knowledge, this is the first real-world application of the Koopman operator-based reduced-order models to a control problem in the field of power electronics. The presented results are very promising, however, it should be stressed out that this contribution is only a first proof of concept since the Koopman matrix was trained offline by using simulation-generated test data from the reference white-box motor model. 

There are plenty of future research questions in terms of using a Koopman operator-based model for MPC. For example, the data collection for Koopman training can be realized using measured motor data from the test bench. By doing so, the Koopman ROM will inherently take real-world effects into account like (cross-)saturation, iron losses effects or flux linkage harmonics. This could lead to improved MPC performance without increasing the computational load due to a more complex plant model. Moreover, the training process and thus the ROM can be adapted online using streaming data \cite{PK18}. Besides that, the motor's torque can be incorporated as an additional observation for the training, and the FCS-MPC can be extended to directly control the torque in an open-loop manner. In a larger scope, the presented modeling method can be also applied to more complex power electronic applications like motor drives including sine-filters between inverter and motor or DC-DC LLC resonant converters where important parasitic effects are hard to model in white-box approaches.

\ifnum\Review=1
	\FloatBarrier
\fi


\bibliographystyle{IEEEtranTIE}
\bibliography{refs}

\end{document}